\documentclass[10pt,draft]{amsart} 
\usepackage{amsmath,amssymb}
\usepackage{amsrefs}

\numberwithin{equation}{section}
\newtheorem{theorem}[equation]{Theorem}
\newtheorem{lemma}[equation]{Lemma}
\newtheorem{corollary}[equation]{Corollary}
\newtheorem{prop}[equation]{Proposition}

\theoremstyle{definition}
\newtheorem{definition}[equation]{Definition}

\theoremstyle{remark}
\newtheorem{remark}[equation]{Remark}

\begin{document}

\begin{center}
\texttt{Comments, suggestions, corrections, and references 
  welcomed!}
\end{center}

\title[Capable groups of prime exponent and class two]{Some results on capable groups of prime exponent and class two}
\author{Arturo Magidin} 
\address{Dept.~Mathematical Sciences, The University of Montana, Missoula MT 59812-0864}
\email{magidin@member.ams.org}


\begin{abstract}
This note collects several results on the capability of $p$-groups of
class two and prime exponent. Among our results, we
settle the $4$-generator case for this class.
\end{abstract}

\maketitle

\section*{Introduction.}

In this note we collect sundry results related to the capability of finite
$p$-groups of class two and prime exponent, using the approach 
introduced in~\cite{capablep}. We restrict to odd primes,
since the case of abelian groups is well understood. I owe a 
considerable amount of the material here to discussions with David
McKinnon\nocite{persdave}, who helped me clarify much of the background
and ideas. 

A group $G$ is said to be \textit{capable} if and only if there exists
a group $K$ such that $G=K/Z(K)$. For $p$-groups (groups of finite
prime-power order) capability is closely related to their
classification. Baer characterized the capable groups which are direct
sums of cyclic groups in~\cite{baer}; the capable extra-special
$p$-groups were characterized by Beyl, Felgner, and Schmid
in~\cite{beyl} (only the dihedral group of order~$8$ and the
extra-special groups of order $p^3$ and exponent~$p$ are capable);
the three authors also describe the metacyclic groups which are capable. The
author characterized the $2$-generated capable $p$-groups of class
two~\cites{capable,twocubed} (for odd~$p$, independently obtained in
part by Bacon and~Kappe in~\cite{baconkappe}).

In~\cite{capablep} we analyzed the special case of $p$-groups of
class two and exponent $p$, where $p$ is an odd prime. We described a way
to translate the problem into a statement of linear algebra,
and obtained several results. We will introduce what I hope is clearer
notation here, reprove some of the results using it, and we will also
be able to obtain new results. Some results from~\cite{capablep} we
will simply quote, and we direct the reader there for any ommitted
proofs.

\section{The set-up.}\label{sec:setup}

We begin with a general construction and some observations.

\begin{definition} Let $V_1$ and $V_2$ be vector spaces,
  and let $\{\ell_i\colon V_1\to V_2\}_{i \in I}$ be a nonempty family of
  linear transformations. 
  For every subspace $X$ of $V_1$, we define ${X}^* < V_2$ to be:
\[ X^* = \Bigl\langle \ell_i(X) \Bigm| i\in I\Bigr\rangle.\]
For every subspace $Y$ of $V_2$, we define $Y^* < V_1$ to be:
\[ Y^* = \bigcap_{i\in I} \ell_i^{-1}(Y).\]
\end{definition}

Note that for all subspaces of $V_1$, if $X\subset X'$, then
$X^*\subset X'^*$; likewise, for subspaces of $V_2$, if $Y\subset
Y'$, then $Y^*\subset Y'^*$.

\begin{theorem} Let $V_1$ and $V_2$ be vector spaces, 
  and let $\{\ell_i\colon V_1\to
  V_2\}_{i\in I}$ be a nonempty family of linear transformations. The operator on the
  subspaces of $V_1$ defined by $X\mapsto X^{**}$ is a closure
  operator; that is, it is increasing, isotone, and idempotent.
  Moreover, if
  $X$ is a subspace of $V_1$, then $X^* = X^{***}$.
\end{theorem}

\begin{proof} It is clear that
  $X\subset X^{**}$, so the operator is increasing; since $X\subset
  X'$ implies $X^*\subset X'^*$, which in turn implies $X^{**}\subset
  X'^{**}$, the operator is isotone.

  To show the operator is idempotent, we want to show that
  $X^{**} = (X^{**})^{**}$. For simplicity, write $X^{**}=Z$. Since
  the operator is increasing, we know that $Z\subset Z^{**}$. By
  construction, we also have that $\ell_i(Z)\subset X^*$ for each $i$,
  so $Z^*\subset X^*$. From this, we get that $Z^{**} \subset
  X^{**}=Z$, as desired.

  Finally, let $X < V_1$. Since $X\subset X^{**}$, we must have $X^*\subset
  X^{***}$. Conversely, from the definition of $X^{***}$ it follows
  that $\ell_i(X^{**})\subset X^*$ for all $i$, so $X^{***}\subset X^*$,
  giving equality.
\end{proof}

The dual property is true for subspaces of $V_2$:

\begin{theorem} Let $V_1$ and $V_2$ be vector spaces, 
  and let $\{\ell_i\colon V_1\to V_2\}_{i\in I}$ be a
  nonempty family of
  linear transformations.
 The operator on subspaces of $V_2$ defined by $Y\mapsto
  Y^{**}$ is an interior operator; that is, it is decreasing, isotone,
  and idempotent. Moreover, if $Y$ is a subspace of $V_2$, then $Y^*=Y^{***}$.
\end{theorem}

\begin{proof} As before, the operator is isotone. Since
  $Y^{**} = \langle \ell_i(Y^*)\rangle$, and 
  $Y^*\subset \ell_i^{-1}(Y)$ for each $i$, it follows that
  $Y^{**}\subset Y$, so the operator is decreasing.

To show the operator is idempotent, set $Z=Y^{**}$. Since $Y^*\subset
Y^{***}$, we must have $Z = Y^{**}\subset Y^{****}=Z^{**}$.
Since the
other inclusion always holds, this shows that $Z=Z^{**}$, as desired.

Finally, since $Y^{**}\subset Y$, it follows that $Y^{***}\subset
Y^*$.  But since the operator ${}^{**}$ on subspaces of $V_1$ is
increasing, we also know that $Y^*\subset Y^{***}$, giving equality.
\end{proof}

\begin{definition} Let $V_1$, $V_2$ be vector spaces, and let
  $\{\ell_i\colon V_1\to V_2\}_{i\in I}$ be a family of linear
  transformations. We will say that a subspace $X$ of $V_1$ is
  \textit{$\{\varphi_i\}_{i\in I}$-closed} (or simply \textit{closed} if there is no
  danger of ambiguity) if and only
  if $X=X^{**}$. 
\end{definition}

Recall that if $V$ is a vector space and $k$ a positive
  integer,  the Grassmannian $Gr(k,V)$ is defined to be the set of all
  $k$-dimensional subspaces of~$V$. It has more algebraic structure
  than simply being a set, but for now we will simply use it as
  convenient notation.

We fix an odd prime $p$ throughout. Let $n$ be an integer
greater than~$1$.  We define three vectors spaces with distinguished bases:

\begin{definition}
Let $U(n)$ be the vector space over the finite field
$\mathbb{F}_p$ of $p$ elements, with basis vectors $x_1,\ldots,x_n$.
Let $V(n)$ be the vector space over $\mathbb{F}_p$ of dimension
$\binom{n}{2}$, with basis vectors $v_{ji}$, $1\leq i<j\leq n$. Let
$W(n)$ be the vector space over $\mathbb{F}_p$ of dimension
$2\binom{n}{2} + 2\binom{n}{3}$, with basis vectors $w_{jik}$, $1\leq
i<j\leq n$, $i\leq k\leq n$. 
\end{definition}

If there is no danger of ambiguity, and $n$ is understood from
context, we will simply use $U$, $V$, and~$W$ instead of $U(n)$,
$V(n)$, and $W(n)$. This will be the case in most of our applications.

We define two families of $n$ linear operators:

\begin{definition}
The linear maps $\psi_i\colon U\to V$, $i=1,\ldots,n$ are defined to be:
\[ \psi_i(x_j) = \left\{\begin{array}{ll}
v_{ji}&\mbox{if $i<j$}\\
-v_{ij}&\mbox{if $i>j$}\\
\mathbf{0}&\mbox{if $i=j$.}
\end{array}\right.\]
\end{definition}

\begin{definition}
The linear maps $\varphi_k\colon V \to W$, $k=1,\ldots,n$,
are defined to be:
\[ \varphi_k(v_{ji}) = \left\{\begin{array}{ll}
w_{jik} & \mbox{if $k\geq i$,}\\
w_{jki}-w_{ikj} & \mbox{if $k<i$.}
\end{array}\right.\]
\end{definition}

Note that there is a natural identification of $V$ with $U\wedge U$,
and of $W$ with $(U\wedge U) \otimes U$ modulo the Jacobi identity:
\[ (u_i\wedge u_j) \otimes u_k + (u_j\wedge u_k)\otimes u_i +
(u_k\wedge u_i)\otimes u_j = 0.\]
Under these identifications, the maps $\psi_i$ correspond to
$\psi_i(\mathbf{u}) = \mathbf{u}\wedge u_i$, and the maps $\varphi_k$
correspond to $\varphi_k(\mathbf{v}\wedge \mathbf{w}) =
(\mathbf{v}\wedge \mathbf{w})\otimes u_k$. 

The definitions of $U$, $V$, $W$, $\psi_j$, and $\varphi_k$ are 
derived from a group theoretic construction. Namely: let $F$ be the
free product of $n$ cyclic groups of order~$p$, generated by
$x_1,\ldots,x_n$ respectively. Let $F_k$ be the $k$-th term of the
lower central series of $F$. 
Then $U$ corresponds to~$F^{\rm ab}=F/F_2$;
$V$ corresponds to $F_2/F_3$, by identifying $[x_j,x_i]$ with
$v_{ji}$; and $W$ corresponds to $F_3/F_4$, identifying
$[[x_j,x_i],x_k]$ with $w_{jik}$. The linear maps $\psi_i$
and~$\varphi_i$ correspond to
the maps induced by the commutator bracket $[-,x_j]$ from $F/F_2$ to
$F_2/F_3$, and from $F_2/F_3$ to $F_3/F_4$, respectively.
We could collect all of
this information by making $U\oplus V\oplus W$ into a Lie algebra, via
the Lie product given by the commutator bracket, but we will not go
through the details here since we will not use the Lie algebra structure.

The connection to capability is explored and explained
 in~\cite{capablep}; we give a very brief summary here:
let $G$ be any
$p$-group of class two and exponent~$p$, with $G^{\rm ab}$ of
rank~$n$. Let $g_1,\ldots,g_n$ be elements of $G$ that project onto a
basis of~$G^{\rm ab}$, and let $K=F/F_3$. Then $K$ is the relatively
free group of rank~$n$ in the variety of nilpotent groups of class at
most two and exponent~$p$, so there is a unique homomorphism $K\to G$
mapping $x_i$ to $g_i$. Let $N$ be the kernel. Then $N\subset
[K,K]$. Thus, $N$ corresponds to a subspace $X$ of~$V$ (since
$[K,K]=F_2/F_3$ corresponds to~$V$). We then have:

\begin{theorem}[Theorem~3.6 in~\cite{capablep}] Fir $n>1$, and let
  $G$, $K$, and $N$ be as in the previous paragraph; identify $[K,K]$
with $V$ by letting $[x_j,x_i]$ correspond to $v_{ji}$, $1\leq i<j\leq
n$, and let $X$ be the subspace of $V$ corresponding to~$N$. Then $G$
is capable if and only if $X$ is closed with respect to the
$\varphi_i$.
\end{theorem}

\section{The linear algebra.}\label{sec:linearalgebra}

We collect here some observations on the spaces $V$ and~$W$, and the
linear transformations $\varphi_k$ defined above.

\begin{definition} Let $i,j$ be integers, $1\leq i< j \leq n$. We
  let $\pi_{ji}\colon V\to \langle v_{ji}\rangle$ be the
  canonical projection. 
\end{definition}

\begin{definition} Let $i,j,k$ be integers, $1\leq i<j\leq n$, $i\leq
  k\leq n$. We let $\pi_{jik}\colon W\to \langle w_{jik}\rangle$ 
  be the canonical projection. 
\end{definition}

\begin{definition} Let $i$ be an integer, $1\leq i\leq n$. We let
  \[\Pi_i\colon V\to \langle v_{i1}, v_{i2},\ldots,
  v_{i(i-1)},v_{(i+1)i},\ldots,v_{ni}\rangle\]
be the canonical projection.
\end{definition}

\begin{lemma} For each $k$, $\varphi_k$ is one-to-one.
\end{lemma}

\begin{lemma} Let $\mathbf{w}\in\varphi_i(V)$. If
  $\pi_{rst}(\mathbf{w})\neq \mathbf{0}$, then $s\leq i\leq t$;
  moreover, exactly one of the inequalities is strict.
\label{lemma:nonzerocoord}
\end{lemma}

\begin{lemma} Fix $n,i,j$, with $n>1$ and $1\leq i<j\leq n$. Then
  $\varphi_i(V)\cap \varphi_j(V) = \{\mathbf{0}\}$.
\end{lemma}

\begin{proof} Assume $\varphi_i(\mathbf{v})\in \varphi_j(V)$, and
  that $\pi_{sr}(\mathbf{v})\neq \mathbf{0}$. If $r\leq i$,
  then $\pi_{sri}(\varphi_i(\mathbf{v}))\neq\mathbf{0}$; but since
  $\varphi_i(\mathbf{v})\in\varphi_j(V)$, this implies that
  $r\leq j\leq i$, which is impossible. 
If, on the other hand, $i<r$, then
$\pi_{sir}(\varphi_i(\mathbf{v}))\neq \mathbf{0}$ and
$\pi_{ris}(\varphi_i(\mathbf{v}))\neq\mathbf{0}$; again, since
$\varphi_i(\mathbf{v})\in\varphi_j(V)$, Lemma~\ref{lemma:nonzerocoord}
implies that $s=r=j$, which is also impossible.
\end{proof}

\begin{lemma} Fix $n>1$, $r\leq n$. Let $i_1,\ldots,i_r$ be integers,
  $1\leq i_1\leq \cdots\leq i_r\leq n$. Then
\[ \varphi_{i_1}(V) \cap \Bigl\langle
\varphi_{i_2}(V),\ldots,\varphi_{i_r}(V)\Bigr\rangle\]
is of dimension $\binom{r-1}{2}$, with basis given by all vectors
$w_{ai_1b}-w_{bi_1a}$, $a,b\in \{i_2,\ldots,i_r\}$, $a> b$. A basis
for the pullback
\[ \varphi_{i_1}^{-1}\left(\Bigl\langle
\varphi_{i_2}(V),\ldots,\varphi_{i_r}(V)\Bigr\rangle\right)\]
is given by $\{v_{ab}\}$, with $a,b \in\{i_2,\ldots,i_r\}$, $a>b$.
\label{lemma:intersections}
\end{lemma}

\begin{proof} The result is true for $r=1$ and $r=2$. Assume then that
  $r\geq 3$. Note that the vectors $w_{ai_1b}-w_{bi_1a}$ are indeed in the
intersection; for any choice of $a>b$ as given, we have 
$\varphi_{i_1}(v_{ab}) = w_{ai_1b} - w_{bi_1a} =
\varphi_b(v_{ai_1})-\varphi_a(v_{bi_1})$.

Let $S = \langle
\varphi_{i_2}(V),\ldots,\varphi_{i_r}(V)\rangle$. 
Assume that $\varphi_{i_1}(\mathbf{v}) \in S$,
and that $\pi_{ts}(\mathbf{v})\neq \mathbf{0}$, where $1\leq s<t\leq n$.
Since $\varphi_{i_1}(\mathbf{v})\in S$, and every vector in~$S$ has
trivial image under all $\pi_{abi_1}$, it follows that $s>i_1$.
Therefore $\varphi_{i_1}(v_{ts}) = w_{ti_1s} - w_{si_1t}$, and so we have:
\[
\varphi_{i_1}(\mathbf{v}) = \varphi_{i_1}\left(\sum_{i_1<s<t\leq n}
\alpha_{ts}v_{ts}\right)
 =  \sum_{i_1<s<t\leq n} \alpha_{ts}(w_{ti_1s}-w_{si_1t}).
\]
Let $\mathbf{w}\in S$; if $\pi_{jik}(\mathbf{w})\neq\mathbf{0}$  with
$i<i_2$, then we must have $k=i_c$ for some $c\in \{2,\ldots,r\}$. Therefore, if
$\varphi_{i_1}(\mathbf{v})\in S$, and $\pi_{ts}(\mathbf{v}) = \alpha_{ts}v_{ts}\neq
\mathbf{0}$, then $s$ and~$t$ are in $\{i_2,\ldots,i_r\}$. Thus we may rewrite
$\varphi_{i_1}(\mathbf{v})$ as:
\[
\varphi_{i_1}(\mathbf{v}) = \sum_{2\leq s<t\leq r}
\alpha_{ts}(w_{i_ti_1i_s}-w_{i_si_1i_t})
 =  \sum_{2\leq s<t\leq r}\varphi_{i_s}(\alpha_{ts}v_{i_ti_1}) -
\varphi_{i_t}(v_{i_si_1}).
\]
This proves the lemma.
\end{proof}

\begin{corollary} Fix $n>1$, $r\leq n$, and let $1\leq i_1\leq
  \cdots\leq i_r \leq n$ be
  integers. Then
\[ \dim\left(\langle
\varphi_{i_1}(V),\ldots,\varphi_{i_r}(V)\rangle\right) = r\binom{n}{2} -
\binom{r}{3}.\]
\end{corollary}

\begin{proof} Each $\varphi_k$ is injective, so we have:
\begin{eqnarray*}
\lefteqn{\dim\left(\langle
\varphi_{i_1}(V),\ldots,\varphi_{i_r}(V)\rangle\right)}\\
 &=&
\left(\sum_{k=1}^r \dim(\varphi_{i_k}(V))\right)
 -
\left(\sum_{k=1}^{r-2}\dim\Bigl(\varphi_{i_{r-k-1}}(V)\cap \langle
\varphi_{i_r}(V),\ldots,\varphi_{i_{r-k}}(V)\rangle \Bigr)\right)\\
& = & r\binom{n}{2} - \left(\sum_{k=1}^r \binom{k-1}{2}\right)
 =  r\binom{n}{2} - \binom{r}{3},
\end{eqnarray*}
as claimed.
\end{proof}

We have concentrated on the ``lowest index'' for simplicity. Of
course, given the definitions, our treatment has symmetry; for example:

\begin{prop} Let $i_1,\ldots,i_r\in\{1,\ldots,n\}$ be pairwise
  distinct. Then
\[ \varphi_{i_1}(V) \cap \Bigl\langle
\varphi_{i_2}(V),\ldots,\varphi_{i_r}(V)\Bigr\rangle\]
has dimension $\binom{r-1}{2}$. Moreover, a basis for the pullback is
given by the vectors $v_{ab}$, with $a,b\in\{i_2,\ldots,i_r\}$, and $a>b$.
\end{prop}

\begin{proof} It is easy to verify that
\[ \varphi_{i_1}(V)\cap \Bigl\langle \varphi_{i_2}(V),\ldots,\varphi_{i_r}(V)\Bigr\rangle\]
is generated by the vectors $w_{abi_1}$ for
$a,b\in\{i_2,\ldots,i_r\}$, $a,i_1>b$, and the vectors
$w_{ai_1b}-w_{bi_1a}$ when $b>i_1$. These vectors are linearly
independent, and pulling them back gives the desired result.
\end{proof}

As we will see below, these descriptions allow us to place bounds on
$\dim(X^*)$ in terms of $\dim(X)$ and $n$.

\section{Dimension counting.}\label{sec:dimcounting}

In this section we will prove some bounds on $\dim(X^*)$ in terms of
$\dim(X)$ and~$n$. These bounds will give sufficient
conditions for $X$ to be closed in terms of its dimension.

\begin{lemma}[Prop.~4.6 in~\cite{capablep}] Fix $n>1$, and 
let $X$ be a subspace of~$V$. If $\dim(X)=1$, then $\dim(X^*)=n$; if
$\dim(X)=2$, then $\dim(X^*)=2n$.
\end{lemma}

In~\cite{capablep} we proved this by considering the collection of
images of a basis of~$X$ under $\varphi_1,\ldots,\varphi_n$, and
showing they would necessarily be linearly independent. We give a new
proof here:

\begin{proof} The results are trivial if $n=2$, so assume $n\geq
  3$. We will prove the contrapositive. Note that 
  $\dim(X^*)\leq n\dim(X)$ always holds, so assume that
$\dim(X^*)<n\dim(X)$. We will show that $\dim(X)\geq 3$.

 Let $k$ be any integer such that
\[ \varphi_k(X) \cap \Bigl\langle
\varphi_{k+1}(X),\ldots,\varphi_n(X)\Bigr\rangle \neq
\{\mathbf{0}\}.\]
Such a $k$ must exist, for if all such intersections were trivial
then we would have $\dim(X^*)=n\dim(X)$. Let 
\[ \mathbf{v}\in \varphi_k^{-1}\Bigl(\bigl\langle
\varphi_{k+1}(X),\ldots,\varphi_{n}(X)\bigr\rangle\Bigr) -
\bigl\{\mathbf{0}\bigr\}.\]
Then from Lemma~\ref{lemma:intersections} and the injectivity of
$\varphi_k$, we deduce that we can write
\[ \mathbf{v} = \!\!\!\!\!\sum_{k+1\leq i<j\leq n} 
\!\!\!\!\!\alpha_{ji}v_{ji},\]
and not all $\alpha_{ji}$ are zero.
Therefore, 
\begin{eqnarray*}
\varphi_k(\mathbf{v}) &=& \sum_{k+1\leq i<j\leq n}
\!\!\!\!\alpha_{ji}(w_{jki}-w_{ikj})\\
& = & \sum_{k+1\leq i<j\leq n} \!\!\!\!\alpha_{ji}(\varphi_i(v_{jk}) -
\varphi_j(v_{ik}))\\
& = & \sum_{i=k+1}^n \varphi_i\left( \sum_{j=i+1}^n \alpha_{ji}v_{jk}
- \sum_{j=k+1}^{i-1} \alpha_{ij}v_{jk}\right).
\end{eqnarray*}
By injectivity of the $\varphi$ and the choice of $\mathbf{v}$,
we again conclude that for each value of $i$, $i=k+1,\ldots,n$,
\[ \mathbf{v}_i = \left(\sum_{j=i+1}^n \alpha_{ji}v_{jk}
- \sum_{j=k+1}^{i-1} \alpha_{ij}v_{jk}\right) \in X.\]
For simplicity, set $\alpha_{ij}=-\alpha_{ji}$, and
$\alpha_{ii}=0$. Then we can rewrite $\mathbf{v}_i$ as:
\[ \mathbf{v}_{i} = \sum_{k+1\leq j\leq n} \!\!\!\!\!\alpha_{ji}v_{jk}.\]
Let $S = \{\mathbf{v}_{k+1},\mathbf{v}_{k+2},\ldots,\mathbf{v}_{n}\}$.
Note that if a subset $T\subset S$ is linearly independent, then so is
$T\cup\{\mathbf{v}\}$.
We will show that that there are
at least two vectors in $S$ that are not multiples of each other;
this will show that the dimension of~$X$ is at least three, proving
the result.

Indeed, fix $i_0$ and $j_0$ such that $\alpha_{j_0i_0}\neq 0$, $k<i_0<j_0\leq n$. 
We have that
\[ \mathbf{v}_{i_0} = \sum_{j=k+1}^n \alpha_{ji_0}v_{jk} \qquad\mbox{and}\qquad
\mathbf{v}_{j_0} = \sum_{j=k+1}^n \alpha_{jj_0}v_{jk},\] so
$\pi_{j_0k}(\mathbf{v}_{i_0})\neq\mathbf{0}$ and
$\pi_{i_0k}(\mathbf{v}_{i_0})=\mathbf{0}$; on the other hand,
$\pi_{i_0k}(\mathbf{v}_{j_0})\neq\mathbf{0}$ and
$\pi_{j_0k}(\mathbf{v}_{j_0})=\mathbf{0}$. So neither of
$\mathbf{v}_{i_0}$ and $\mathbf{v}_{j_0}$ is a scalar multiple of the
other, and $\{\mathbf{v}_{i_0},\mathbf{v}_{j_0},\mathbf{v}\}$ is a
linearly independent subset of~$X$. Therefore, the dimension of~$X$
must be at least $3$, as claimed. \end{proof}

With similar arguments we can obtain a (naive) lower bound for
the dimension of $X^*$ in terms of $\dim(X)$. We have:

\begin{lemma} Fix $n\geq 3$, and let $X$ be a subspace of~$V$, with
  $\dim(X)=k$. Then 
\[ nk \geq \dim(X^*)\geq \left\{\begin{array}{lcl}
\displaystyle nk - \left(\frac{k-1}{2}\right)^2&&\mbox{if $k$ is odd,}\\
\\
\displaystyle nk - \left(\frac{k}{2}\right)\left(\frac{k-2}{2}\right)&&\mbox{if $k$ is even.}
\end{array}\right.\]
\label{lemma:bounds}
\end{lemma}

\begin{proof} The upper bound is trivial. 
Let
\[ Z_i = \varphi_i(X)\cap \Bigl\langle
\varphi_{i+1}(X),\ldots,\varphi_n(X)\Bigr\rangle\] for
$i=1,2,\ldots,n-2$. Note that
\[ \dim(X^*) = nk - \sum_{i=1}^{n-2}\dim(Z_i).\]

Let $1\leq i_1<i_2<\cdots<i_r\leq n-2$ be the indices for which $Z_i$
is nontrivial.  The dimension of $Z_{i_1}$ is at most $k-2$; this
since, as above, if $\mathbf{v}$ is in $\varphi_{i_1}^{-1}(Z_{i_1})$,
then $\Pi_{k}(\mathbf{v})\neq \mathbf{0}$ implies
$k\in\{i_2,\ldots,i_r\}$; 
and we need at least two linearly independent vectors with nontrivial
image under $\Pi_{i_1}$ to account for the intersection. Likewise, the
dimension of $Z_{i_2}$ is at most $k-4$, since the two vectors with
nontrivial image under $\Pi_{i+1}$ cannot be involved in this
intersection. Continuing this way, we obtain
\[ \sum_{j=1}^r \dim(Z_{i_j}) \leq (k-2) + (k-4) + \cdots + (k-2r).\]
Moreover, note that $2r< k$. If $k$ is odd, say $k=2\ell+1$, then this
gives
\begin{eqnarray*}
\sum_{j=1}^r \dim(Z_{i_j}) &\leq& (k-2) + (k-4) + ... + (k-2\ell)\\
& = & \ell k - 2(1+2+\cdots+\ell) \\
& = & 2\ell^2 + \ell - \ell^2 - \ell\\
& = & \left(\frac{k-1}{2}\right)^2.
\end{eqnarray*}
On the other hand, if $k$ is even, $k=2\ell$, then
\begin{eqnarray*}
\sum_{j=1}^r \dim(Z_{i_j}) &\leq& (k-2) + (k-4) + ... + (k-2(\ell-1))\\
& = & (\ell-1)k - 2(1+2+\cdots+(\ell-1)) \\
& =& 2\ell^2 - 2\ell - \ell^2 + \ell\\ 
& = & \left(\frac{k}{2}\right)^2 - \frac{k}{2}\\
& = & \left(\frac{k}{2}\right)\left(\frac{k-2}{2}\right).
\end{eqnarray*}
This proves the lower bound.
\end{proof}

\begin{remark} We can improve upon the lower bound with more
  care; for example, as we will see below,
  whereas the lemma only gives $\dim(X^*)\geq 16$ in the case $n=4$ and
  $k=5$, we can actually prove that we will have $\dim(X^*)\geq 17$.
\end{remark}

Though very rough, the lower bounds established below give us a
sufficient condition for $X$ to be closed in terms only of
$\dim(X)$ and~$n$:

\begin{theorem}
Let $n>1$, and let $X$ be a subspace of~$V$. If $\dim(X)=k$, then $X$
is closed whenever $k<2\sqrt{n}$.
\label{th:sufficientdimension}
\end{theorem}

\begin{proof} 
Let $X'$ be a subspace of dimension $k+1$; we show that
  if $k$ satisfies the given inequality  then
  $\dim(X'^*)>\dim(X^*)$. 

From Lemma~\ref{lemma:bounds} we know that $\dim(X^*)\leq nk$ and 
\[ \dim(X'^*)\geq \left\{\begin{array}{lcl}
\displaystyle n(k+1) - \left(\frac{k+1}{2}\right)\left(\frac{k-1}{2}\right)&&\mbox{if $k$ is odd,}\\
\\
\displaystyle n(k+1) - \left(\frac{k+1}{2}\right)^2&&\mbox{if $k$ is even.}
\end{array}\right.\]
Thus, $\dim(X^*)<\dim(X'^*)$ will certainly hold provided that
\[\begin{array}{rclcl}
k^2 & < & 4n + 1 &&\mbox{if $k$ is odd;}\\
k^2 & < & 4n &&\mbox{if $k$ is even.}
\end{array}\]
It is easy to verify that these inequalities hold if and only if
$k<2\sqrt{n}$, regardless of the parity of~$k$.
This proves that $\dim(X^*)<\dim(X'^*)$ for any $X'$ of dimension
$k+1$, provided $k<2\sqrt{n}$.

Since any $X'$ that satisfies $X\subseteq X'\subseteq X^{**}$ must
have $X'^*=X^*$, it follows that there are no subspaces of dimension
$k+1$ that lie between $X$ and $X^{**}$; this means that $X^{**}=X$,
so $X$ is closed, as claimed.
\end{proof}

Translating into the group theoretic setting, we have:

\begin{theorem}
Let $G$ be a $p$-group of class two and exponent~$p$, ${\rm
  rank}(G^{\rm ab})=n$, and ${\rm rank}([G,G])=\ell$.
If $\ell > \binom{n}{2} - 2\sqrt{n}$ then $G$ is capable.
\label{th:largeenoughcomm}
\end{theorem}

We also have:

\begin{corollary}[cf.~Cor.~3.7 in~\cite{capablep}]
Let $n\geq 3$, and $X$ be a subspace of~$V$. If $\dim(X)\leq 4$, then
$X$ is closed.
\label{cor:uptofour}
\end{corollary}

\begin{proof} Let $\dim(X)=k$. Theorem~\ref{th:sufficientdimension}
  gives the result if $k\leq 3$, or if $k=4<n$. So we may assume
  $\dim(X)=n=4$.
From Lemma~\ref{lemma:bounds}, we know that if $\dim(X')=5$, then
\[16  \leq  \dim(X'^*)  \leq  20.\]
We claim that in fact we can sharpen the lower bound to
$\dim(X'^*)\geq 17$. If so this will prove the proposition via the
same argument as in the final paragraph of the proof of
Theorem~\ref{th:sufficientdimension}.  There are only two possible
values of $i$ for which
\[Z'_i = \varphi_i(X')\cap \Bigl\langle
\varphi_{i+1}(X'),\ldots,\varphi_n(X')\Bigr\rangle\]
will be nontrivial: $i=1$ and $i=2$. The dimension
of $Z'_2$ is at most $1$, and the dimension of $Z'_1$ is at most
$3$; but the only way for $Z'_1$ to be of dimension $3$ is if $X'$
contains all of $v_{32}$, $v_{42}$, and $v_{43}$, in which case it is
trivial to see that we must also have $v_{21}$, $v_{31}$, and $v_{41}$
in~$X'$. Since $\dim(X')=5$, this is impossible. Thus, in
fact, $\dim(Z'_2)\leq 1$ and $\dim(Z'_1)\leq 2$, so $\dim(X'^{*})\geq
20-3 = 17$. This proves the claim, and the proposition.
\end{proof}

\section{Some results.}\label{sec:oldresults}

We quote the following two results without proof:

\begin{theorem}[Lemma 5.3 in~\cite{capablep}] If $X$ is a coordinate
  subspace (that is, it is generated by a subset of the $v_{ji}$) then
  $X$ is closed.
\label{th:coordinate}
\end{theorem}

\begin{theorem}[Lemma 5.7 in~\cite{capablep}]
Let $m$ be fixed, $1<m<n$, and assume that $X_1$ is a
  subspace of $\langle v_{ji}\,|\, 1\leq i<j\leq m\rangle$, and $X_2$
  is a subspace of $\langle v_{ji}\,|\, m<i<j\leq n\rangle$. Then
  $X_1\oplus X_2$ is closed.
\label{th:smallplusbig}
\end{theorem}

In~\cite{capablep}*{Lemma~5.9} we proved a special case of the following result:

\begin{lemma}
Let $n>1$ be an integer, let $I$ be a proper nonempty subset of
$\{1,\ldots,n\}$, and let $J= \{1,\ldots,n\}-I$. Let 
\begin{eqnarray*}
V_I &=& \langle v_{ji}\,|\, 1\leq i<j\leq n;\ i,j\in I\}\\
V_J &=& \langle v_{ji}\,|\, 1\leq i<j\leq n;\ i,j\in J\}\\
V_{(I,J)} & = & \langle v_{ji}\,|\,1\leq i<j\leq n;\ \mbox{exactly one
  of $i,j$ is in $I$, one in $J$}\rangle.
\end{eqnarray*}
Let $X_I$ be a subspace of $V_I$, $X_J$ be a subspace of
$V_J$, and let
\[ X = X_I \oplus X_J \oplus V_{(I,J)}.\]
Then 
\[ X^{**} = \left(\bigcap_{i\in
  I}\varphi_i^{-1}\Bigl(\bigl\langle \varphi_k(X_I)\,\bigm|\,k\in I\bigr\rangle\Bigr)\right)
    \oplus \left(\bigcap_{j\in
  J}\varphi_j^{-1}\Bigl(\bigl\langle \varphi_k(X_J)\,\bigm|\, k \in
  J\bigr\rangle\Bigr)\right)\oplus V_{(I,J)}.\]
In particular, $X$ is
  $\{\varphi_k\}_{k=1}^n$-closed if and only if $X_I$ is
  $\{\varphi_i\}_{i\in I}$-closed and $X_J$ is $\{\varphi_j\}_{j\in
  J}$-closed.
\end{lemma}

\begin{proof} Let $W'$ be the subspace of $W$ generated by all $w_{jik}$ in which
exactly one or two of $i,j,k$ are in~$I$. 
First, we claim that $W'\subset X^*$.
Suppose that exactly one of $i,j,k$ lies in~$I$. If $j\in I$, then
$v_{ji} \in V_{(I,J)}$, so $\varphi_k(v_{ji})=w_{jik}\in X^*$. If
$i\in I$, then again $v_{ji}\in V_{(I,J)}$, so again we conclude that
$w_{jik}\in X^*$. Finally, if $k\in I$, then either $v_{jk}$ or $v_{kj}$ lies
in $V_{(I,J)}$ (whichever is appropriate); the image under $\varphi_i$
of this vector is $\pm (w_{jik}-w_{kij})$. Since $w_{kij}\in
X^*$ by the argument above it follows
that $w_{jik}\in X^*$ as desired.
The case when exactly two if $i,j,k$ lie in $I$ follows by symmetry,
since exactly one of them will lie in~$J$.

Next, we claim that
\[ X^* = \Bigl\langle \varphi_i(X_I)\,\Bigm|\, i\in I\Bigr\rangle
\oplus \Bigl\langle \varphi_j(X_J)\,\Bigm|\, j\in J\Bigr\rangle \oplus
W'.\]
Indeed, the right hand side is contained in $X^*$.
It suffices to show that $\varphi_j(X_I)\subset W'$ for every $j\in J$
(the symmetrical argument will show it for $\varphi_i(X_J)$, $i\in I$).
Since $\pi_{abc}\left(\varphi_j(X_I)\right)\neq \{\mathbf{0}\}$ 
if and only if exactly two of $a,b,c$ are in~$I$ and either $b$ or $c$
are equal to~$j$, this establishes our second claim.

Finally, we establish the equality in the lemma: let
\[ \mathbf{v} \in \left(\bigcap_{i\in
  I}\varphi_i^{-1}\Bigl(\bigl\langle \varphi_k(X_I)\,\bigm|\,k\in I\bigr\rangle\Bigr)\right)
     \oplus \left(\bigcap_{j\in
  J}\varphi_j^{-1}\Bigl(\bigl\langle \varphi_k(X_J)\,\bigm|\, k \in
  J\bigr\rangle\Bigr)\right)\oplus V_{(I,J)}.\]
Then we can write $\mathbf{v} = \mathbf{x}_I \oplus \mathbf{x}_{(I,J)}
  \oplus \mathbf{x}_J$, where
\begin{eqnarray*}
\mathbf{x}_I & \in & \bigcap_{i\in I}\varphi_i^{-1}\Bigl(\bigl\langle
\varphi_k(X_I)\,|\,k\in I\bigr\rangle\Bigr),\\
\mathbf{x}_J & \in & \bigcap_{j\in J}\varphi_j^{-1}\Bigl(\bigl\langle
\varphi_k(X_J)\,|\, k\in J\bigr\rangle\Bigr),\\
\mathbf{x}_{(I,J)} & \in & V_{(I,J)}.
\end{eqnarray*}
Since $V_{(I,J)}\subset X \subset X^{**}$, to show that $\mathbf{v}\in
X^{**}$ it is enough to show that $\mathbf{x}_I$ and $\mathbf{x}_J$
are both in~$X^{**}$. And indeed: notice that
$\pi_{ab}(\mathbf{x}_I)\neq\mathbf{0}$ only if $a,b\in I$. Therefore,
$\varphi_j(\mathbf{x}_I)\in W'$ for every $j\in J$. And by
construction we have $\varphi_i(\mathbf{x}_I)\in X^*$ for every $i\in
I$. Thus $\varphi_k(\mathbf{x}_I)\in X^*$ for all~$k$. Symmetrically,
$\varphi_k(\mathbf{x}_J)\in X^*$ for all~$k$ as well. Therefore, each
of $\mathbf{x}_I$ and $\mathbf{x}_J$ lie in $X^{**}$, as desired.

For the converse inclusion, let $\mathbf{x}\in X^{**}$. We can write
$\mathbf{x}=\mathbf{v}_I \oplus\mathbf{v}_J\oplus \mathbf{v}_{(I,J)}$,
with $\mathbf{v}_I\in V_I$, $\mathbf{v}_J\in V_J$, and
$\mathbf{v}_{(I,J)}\in V_{(I,J)}$.  We must have
$\varphi_k(\mathbf{v})\in X^*$ for every $k$. If $k\in I$, then again
we must have
\[\varphi_k(\mathbf{v}_I) \in \Bigl\langle \varphi_i(X_I)\,\Bigm|\,
i\in I\Bigr\rangle,\qquad \varphi_k(\mathbf{v}_J),\varphi_k(\mathbf{v}_{(I,J)})\in W'.\]
And if $k\in J$, then
\[ \varphi_k(\mathbf{v}_J)\in \Bigl\langle \varphi_j(X_J)\,|\Bigm|\,
j\in J\Bigr\rangle,\qquad
\varphi_k(\mathbf{v}_I),\varphi_k(\mathbf{v}_{(I,J)})\in W'.\]
Therefore, 
\begin{eqnarray*}
\mathbf{v}_I & \in & \bigcap_{i\in I}\varphi_i^{-1}\Bigl(\bigl\langle
\varphi_k(X_I)\,|\,k\in I\bigr\rangle\Bigr),\\
\mathbf{v}_J & \in & \bigcap_{j\in J}\varphi_j^{-1}\Bigl(\bigl\langle
\varphi_k(X_J)\,|\, k\in J\bigr\rangle\Bigr),
\end{eqnarray*}
as desired. This proves the equality.

The final clause of the lemma follows from the observation that
$X$ will be closed if and only if 
\[X_I = \bigcap_{i\in I} \varphi_i^{-1}\Bigl(\bigl\langle
\varphi_k(X_I)\,\bigm|\, k\in I\bigr\rangle\Bigr)\ \mbox{and}\ X_J =
\bigcap_{j\in J}\varphi_j^{-1}\Bigl(\bigl\langle
\varphi_k(X_J)\,\bigm|\, k\in J\bigr\rangle\Bigr).\]
\end{proof}

Notice that the subspace $V_{(I,J)}$, when interpreted as a subgroup
of commutators we are making trivial, says that each $x_i$ commutes
with each $x_j$, where $i\in I$ and $j\in J$. From this we obtain the following:

\begin{corollary} Let $G_1$ and $G_2$ be two noncyclic $p$-groups of
  class at most two and exponent~$p$. Then $G=G_1\oplus G_2$ is
  capable if and only if each of $G_1$ and~$G_2$ are capable.
\end{corollary}

The corollary is clearly not true if we drop the ``noncyclic''
hypothesis, since a cyclic group of order~$p$ is not capable, but the
direct sum of two cyclic groups of order~$p$ is capable. To understand
why we must make the distinction, note that when $|I|=1$, $V_I$ is
trivial so $X_I$ is perforce $\{\varphi_i\}_{i\in I}$-closed; but our
developement always assumes $n>1$. The symmetric situation happens
when $|I|=n-1$. However, what we obtain in this instance is the
following result, that is perhaps one of the most unexpected results
from~\cite{capablep}: 

\begin{theorem}[Lemma~5.9 in~\cite{capablep}] Let $n>2$, and let $X$ be a subspace of $V(n-1)$
  and let $X' = X\oplus \langle
  v_{ni}\,|\, 1\leq i\leq n-1\rangle$ be a subspace of $V(n)$. Then
  $X'$ is $\{\varphi_i\}_{i=1}^{n}$-closed if and only if $X$ is
  $\{\varphi_i\}_{i=1}^{n-1}$-closed.
\label{th:cancelcentral}
\end{theorem}

The usefulness of this theorem becomes apparent when we translate it
back into group theoretic terms:

\begin{corollary}  Let $K$ be a finite noncyclic group of exponent $p$ and
  class $2$, and let $G=K\oplus C_p$, where $C_p$ is the
  cyclic group of order~$p$. Then $G$ is capable if and only if $K$ is
  capable.
\end{corollary}

Since every finite group of class at most two and exponent $p$ may be
written as $K\oplus C_p^r$ for some $r\geq 0$ and $K$ satisfying
$[K,K]=Z(K)$, we conclude that:

\begin{corollary} Let $G$ be a $p$-group of class at most two and
  exponent~$p$. Write $G=K\oplus C_p^r$, where $K$ satisfies
  $[K,K]=Z(K)$. Then $G$ is capable if and only if (i) $K$ is nontrivial
  and capable, or (ii) $K$ is trivial and $r\geq 2$.
\end{corollary}

This allows us to ignore the ``extra'' elements in the center of~$G$
that do not result from commutators. Thus, for example, we can
strengthen Theorem~\ref{th:largeenoughcomm} by replacing the rank of
$G^{\rm ab}$ by the rank of $G/Z(G)$. We then have the following:

\begin{corollary}
Let $G$ be a $p$-group of class two and exponent~$p$, ${\rm
  rank}(G/Z(G))=n$, and ${\rm rank}([G,G])=\ell$.
If $\ell > \binom{n}{2} - 2\sqrt{n}$, then 
then $G$ is capable.
\label{cor:suffcardcond}
\end{corollary}

In~\cite{heinnikolova}, the authors established a minimum size for the
commutator subgroup of a capable group~$G$ of class two and
exponent~$p$ which satisfies $[G,G]=Z(G)$; namely, they proved:

\begin{theorem}[Theorem~1 in~\cite{heinnikolova}] Suppose that $G$ is a group of
  exponent~$p$ which satisfies $[G,G]=Z(G)$. If $G$ is capable and
  $[G,G]$ is of rank~$\ell$, then $G/Z(G)$ is of rank at most $2\ell +
  \binom{\ell}{2}$.
\label{th:heinnik}
\end{theorem}

First, our development above now shows that if $G$ is not cyclic, then
Theorem~\ref{th:heinnik} remains true if we drop the hypothesis that
$[G,G]=Z(G)$. Second, note that Theorem~\ref{th:heinnik} gives a
\textit{necessary} condition. If we let
$n = {\rm rank}(G/Z(G))$, as we did in
Corollary~\ref{cor:suffcardcond}, then Theorem~\ref{th:heinnik} gives
the inequality
\[ \ell \geq \frac{\sqrt{9 + 8n} - 3}{2}\]
as a necessary condition for the capability of~$G$.
Combining it with
Corollary~\ref{cor:suffcardcond}, we have:

\begin{theorem} Let $G$ be a noncylic group of exponent~$p$ and class
  at most two. Let ${\rm rank}(G/Z(G))=n$ and ${\rm rank}([G,G])=\ell$.
A necessary condition for the capablity of~$G$ is that
\[ \ell \geq \frac{\sqrt{9+8n}-3}{2}.\]
A sufficient condition for the capability of~$G$ is that
\[ \ell > \binom{n}{2} - 2\sqrt{n}.\]
\end{theorem}

Although Theorem~\ref{th:cancelcentral} focuses on the particular
subspace that results when we take $I=\{1,2,\ldots,n-1\}$, in which
chase $V_{(j,i)}$ is none other than $\langle x_n\rangle^*$ (relative to
$\{\psi_i\}_{i=1}^{n}$), there is nothing special about $x_n$ (except
that the theorem is easy to state). Because of the symmetry inherent
in $U$, we can replace $\langle x_n\rangle^*$ with any other subspace
of the form $\langle\mathbf{u}\rangle^*$, with
$\langle\mathbf{u}\rangle\in Gr(1,U)$. It is easy to verify that in
this case we will have $\langle \mathbf{u}\rangle\in Gr(n-1,V)$. Via a
change of basis for~$U$ and the identifications $V= U \wedge U$,
etc. noted at the end of Section~\ref{sec:setup}, we obtain:

\begin{theorem} Fix $n>2$, and let $X'$ be a subspace of $V(n)$. 
 If there exists $\langle\mathbf{u}\rangle\in Gr(1,U)$ such that
  $\langle\mathbf{u}\rangle^*\subset X'$, then there exists a subspace
  $X$ of $V(n-1)$ of dimension
  $\dim(X')-n+1$ such that $X'$ is $\{\varphi_i\}_{i=1}^{n}$-closed
  if and only if $X$ is $\{\varphi_i\}_{i=1}^{n-1}$-closed.
\label{th:cancelanycentral}
\end{theorem}

\section{The case $n=4$.}\label{sec:nequalfour}

In this section we will settle the 
  case $n\leq 4$.  First, notice that for $n\leq 3$,
  Corollary~\ref{cor:uptofour} guarantees that all subspaces $X$ will
  be closed. For $n=4$, the same corollary settles all cases except
  $\dim(X)=5$ and $\dim(X)=6$; the latter case is just $X=V$, in which
  case $X^{**}=V=X$, so $X$ is closed.
Thus, we may restrict ourselves to $\dim(X)=5$. Then the only two
possibilities are $X^{**}=X$ and $X^{**}=V$. Since $X^{**}=V$ if and
only if $X^*=W$, we have:

\begin{prop}
Let $n=4$, and let $X$ be a $5$-dimensional subspace of~$V$. Then $X$
is closed if and only if $X^*$ is a proper subspace of~$W$.
\end{prop}

Since $\dim(W)=20$, and in this case we have $\dim(W)=n\dim(V)$, the following observation
becomes relevant: 

\begin{theorem}
Let $n>1$ be an integer. Let $X$ is a subspace of $V$, and assume that $\dim(X)\geq n$
 and $\langle \mathbf{u}\rangle^*\subset X$ for some
$\langle\mathbf{u}\rangle \in Gr(1,U)$. Then $\dim(X^*) < n\dim(X)$.
\label{th:dropifwecancancel}
\end{theorem}

\begin{proof}
By symmetry, it is enough to establish this result when
$\mathbf{u}=x_1$. Then $\langle\mathbf{u}\rangle^*$ is generated by
the vectors $v_{21}$, $v_{31},\ldots,v_{n1}$. 

Let $\mathbf{v}\neq \mathbf{0}$ be a vector in $X$ with $\Pi_1(\mathbf{v})=\mathbf{0}$.
Writing
\[ \mathbf{v} = \sum_{2\leq i<j\leq n}\!\!\!\!\! \alpha_{ji}v_{ji}\]
we have:
\begin{eqnarray*}
\varphi_1(\mathbf{v}) & = & \sum_{2\leq i<j\leq
  n}\alpha_{ji}\left(v_{j1i}-v_{i1j}\right)\\
& = & \sum_{2\leq i<j\leq n}\left(\varphi_i(\alpha_{ji}v_{j1}) -
  \varphi_j(\alpha_{ji}v_{i1})\right),
\end{eqnarray*}
so $\varphi_1(\mathbf{v})\in\langle
\varphi_2(X),\ldots,\varphi_n(X)\rangle$. This shows that
$\dim(X^*)<n\dim(X)$, as claimed.
\end{proof}

For $n=4$ and $\dim(X)=5$, we claim that the converse holds.

\begin{theorem}
Let $n=4$, and let $X$ be a $5$-dimensional subspace of~$V$. Then
$\dim(X^*)<20$ if and only if there exists a nonzero vector $\mathbf{u}\in U$
such that $\langle \mathbf{u}\rangle^*\subset X$.
\label{th:dimfivecase}
\end{theorem}

\begin{proof} The if clause follows from
  Theorem~\ref{th:dropifwecancancel}. What we need to do is prove the
  converse. So assume that $\dim(X^*)<20$. Let
\[ Z_1 = \varphi_1(X)\cap\Bigl\langle
  \varphi_2(X),\varphi_3(X),\varphi_4(X)\Bigr\rangle\quad\mbox{and}\quad
  Z_2 = \varphi_2(X)\cap\Bigl\langle
  \varphi_3(X),\varphi_4(X)\Bigr\rangle.\]  
At least one of $Z_1$ and $Z_2$ must be nontrivial.

If $Z_2$ is nontrivial, then it is dimension~$1$, and the only
possibility is for $v_{43}$, $v_{42}$, and $v_{32}$ to be elements
of~$X$.
Let $\mathbf{v}=\alpha_{21}v_{21} + \alpha_{31}v_{31} +
\alpha_{41}v_{41}\in X$ be a nonzero element. This is possible since
$\dim(X)=5$. We claim that $\langle\mathbf{u}\rangle^*\subset X$,
where $\mathbf{u} = \alpha_{21}x_2 + \alpha_{31}x_3 +
\alpha_{41}x_4$. Indeed, we have:
\begin{eqnarray*}
\psi_1(\mathbf{u}) & = & \alpha_{21}v_{21} + \alpha_{31}v_{31} +
\alpha_{41}v_{41} = \mathbf{v} \in X;\\
\psi_2(\mathbf{u}) & = & \alpha_{31}v_{32} + \alpha_{41}v_{42} \in X;\\
\psi_3(\mathbf{u}) & = & -\alpha_{21}v_{32} + \alpha_{41}v_{43} \in X;\\
\psi_{4}(\mathbf{u}) & = & -\alpha_{21}v_{42} - \alpha_{31}v_{43} \in X.
\end{eqnarray*}
Since each $\psi_i(\mathbf{u})\in X$, it follows that
$\langle\mathbf{u}\rangle^*\subset X$, as claimed.

If $Z_1$ is nontrivial, then there must be at least one nonzero vector in
\[\langle v_{32}, v_{42}, v_{43}\rangle\cap X\]
whose image under
$\varphi_1$ lies in
$Z_1$. Say $\mathbf{v}_1= \alpha v_{32} + \beta v_{42} + \gamma v_{43}$,
with $\alpha,\beta,\gamma$ not all zero.  Then $X$ must also contain
the following three vectors:
\begin{eqnarray*}
\mathbf{v}_2 &  = & \alpha v_{31} + \beta v_{41},\\
\mathbf{v}_3 & = & \alpha v_{21}  - \gamma v_{41},\\
\mathbf{v}_4 & = & \beta v_{21} + \gamma v_{31},
\end{eqnarray*}
in order to obtain $\varphi_1(\mathbf{v}_1) = \varphi_2(\mathbf{v}_2) -
\varphi_3(\mathbf{v}_3) - \varphi_4(\mathbf{v}_4)$.

However, $\mathbf{v}_2$, $\mathbf{v}_3$, and $\mathbf{v}_4$ are not
linearly independent. For example, if $\gamma\neq 0$, then 
$\alpha\mathbf{v}_4 + \beta\mathbf{v}_3 = \alpha\gamma v_{31} -
\beta\gamma v_{41} = \gamma(\alpha v_{31} +\beta v_{41}) =
\gamma\mathbf{v}_2$,
so $\mathbf{v}_2$ lies in
$\langle\mathbf{v}_3,\mathbf{v}_4\rangle$. Likewise, if $\beta\neq 0$
then $\mathbf{v}_3$ is a linear combination of $\mathbf{v}_2$ and
$\mathbf{v}_4$; and if $\alpha\neq 0$ then $\mathbf{v}_4$ is a linear
combination of $\mathbf{v}_2$ and $\mathbf{v}_3$. 

If $\alpha\neq 0$, then
taking a suitable scalar multiple of $\mathbf{v}$ we may assume
$\alpha=1$. So $X$ contains the following three linearly independent
vectors:
\begin{eqnarray*}
\mathbf{v}_1 & = & v_{32} + \beta v_{42} + \gamma v_{43},\\
\mathbf{v}_2 & = & v_{31} + \beta v_{41},\\
\mathbf{v}_3 & = & v_{21} - \gamma v_{41}, 
\end{eqnarray*}
for some $\beta,\gamma$ (possibly equal to zero). Since $\dim(X)=5$,
we must have at least two other linearly independent vectors. Adding
suitable multiples of the $\mathbf{v}_j$, $j=1,2,3$, we may
assume that they both lie in $\langle v_{41},v_{42},v_{43}\rangle$.

If $v_{41}\in X$, then $X$ contains $v_{21}$, $v_{31}$, and
$v_{41}$, which means 
\[\langle x_1\rangle^* = \langle
v_{21},v_{31},v_{41}\rangle\subset X.\]
If $v_{42},v_{43}\in X$, then
$v_{32}\in X$ as well, so $Z_2$ is nontrivial and as seen above there
is some $\langle\mathbf{u}\rangle^*\subset X$. So we may assume that
$X$ has a basis consisting of $\mathbf{v}_1$, $\mathbf{v}_2$, and
$\mathbf{v}_3$, together with
two other vectors of the form:
\begin{eqnarray*}
\mathbf{v}_5 & = & v_{41} + \delta v_{43},\\
\mathbf{v}_6 & = & v_{42} + \lambda v_{43}.
\end{eqnarray*}
In this case, let $\mathbf{u} = x_2 + \lambda x_3
+ (\beta\lambda-\gamma) x_4$. Then we have:
\begin{eqnarray*}
\psi_1(\mathbf{u})& = & v_{21} + \lambda v_{31} +
(\beta\lambda-\gamma)v_{41}
 =  \mathbf{v}_3 + \lambda\mathbf{v}_2;\\
\psi_2(\mathbf{u}) & = & \lambda v_{32} +
(\beta\lambda-\gamma)v_{42}
 =  \lambda\mathbf{v}_1 - \gamma\mathbf{v}_{5};\\
\psi_3(\mathbf{u}) & = & -v_{32} + (\beta\lambda-\gamma)v_{43}
 =  -\mathbf{v}_1 + \beta\mathbf{v}_5;\\
\psi_4(\mathbf{u}) & = & -v_{42} -\lambda v_{43} = -\mathbf{v}_5.
\end{eqnarray*}
Thus, $\langle\mathbf{u}\rangle^*\subset X$.

If $\alpha=0$ and $\beta\neq 0$, then we may assume $\beta=1$, so now
the following vectors lie in~$X$:
\begin{eqnarray*}
\mathbf{v}_1 & = & v_{42} + \gamma v_{43},\\
\mathbf{v}_2 & = & v_{41},\\
\mathbf{v}_4 & = & v_{21}+\gamma v_{31}.
\end{eqnarray*}
Two other vectors, which we may assume lie in $\langle
v_{31},v_{32},v_{43}\rangle$, will form a basis for~$X$. If $v_{43}\in
X$, then $\langle x_4\rangle^* = \langle v_{31},v_{32},v_{43}\rangle\subset
X$,
and we are done. Otherwise, we may assume that 
\begin{eqnarray*}
\mathbf{v}_5 & = & v_{31} + \delta v_{43},\\
\mathbf{v}_6 & = & v_{32} + \lambda v_{43},
\end{eqnarray*}
lie in $X$, for some $\delta,\lambda$. In this case, let
$\mathbf{u}=x_2 + \gamma x_3 - \lambda x_4$, and we obtain
\begin{eqnarray*}
\psi_1(\mathbf{u}) & = & \mathbf{v}_4 - \lambda\mathbf{v}_2,\\
\psi_2(\mathbf{u}) &= & \gamma \mathbf{v}_6 - \lambda\mathbf{v}_1,\\
\psi_3(\mathbf{u}) & = & -\mathbf{v}_6,\\
\psi_4(\mathbf{u}) & = & -\mathbf{v}_1.
\end{eqnarray*}
Finally, if $\alpha=\beta=0$, we may assume $\gamma=1$, so now we have
the following linearly independent vectors in~$X$:
\begin{eqnarray*}
\mathbf{v}_1 & = & v_{43},\\
\mathbf{v}_3 & = & -v_{41},\\
\mathbf{v}_4 & = & v_{31},
\end{eqnarray*}
plus two more vectors to form a basis, these latter vectors taken from
$\langle v_{21}, v_{32}, v_{42}\rangle$. If $v_{21}\in X$, then 
$\langle x_1\rangle^* = \langle v_{21}, v_{31},
v_{41}\rangle\subset X$,
and again we are done. Otherwise, we may take
\begin{eqnarray*}
\mathbf{v}_5 & = & v_{32} + \delta v_{21},\\
\mathbf{v}_6 & = & v_{42} + \lambda v_{21},
\end{eqnarray*}
as two other vectors in $X$. This time, letting $\mathbf{u}=-\delta
x_1 + x_3$ gives $\langle \mathbf{u}\rangle^* \subset X$, as desired.
\end{proof}

\begin{corollary} Let $n=4$, and let $X$ be a $5$-dimensional subspace of $V$.
  Then $X$ is closed if and only if there exists a
  nonzero vector $\mathbf{u}\in U$ such that $\langle
  \mathbf{u}\rangle^*\subset X$.
\label{cor:dimxfivewithcentral}
\end{corollary}

This result may seem a bit artificial when put in this guise. However,
translating it back to group theory reveals what is going on. Recall
that the vectors in $U$ correspond to elements of $G^{\rm ab}$, and if
$\langle\mathbf{u}\rangle^*\subset X$, that means that the
corresponding element of $G$ is central (since commutators are central
in~$G$, the ambiguity in the statement ``the corresponding element
of~$G$'' is immaterial). So, in group theoretic language,
Corollary~\ref{cor:dimxfivewithcentral} states:

\begin{theorem} Let $G$ be a group of class two and exponent
  $p$, $p$ and odd prime. If $G$ is of order $p^5$ and $[G,G]$ is
  cyclic of order~$p$, then $G$ is capable if and only if $[G,G]\neq
  Z(G)$. That is, in this case, $G$ is capable if and only if $G$ is not extra-special.
\end{theorem}

Recall that we say a group $G$ is \textit{$n$-generated} to mean
  that it can be generated by $n$ elements, though it may in fact need
  fewer. We say $G$ is \textit{minimally $n$-generated} to mean that
  $G$ can be generated by $n$ elements, but cannot be generated by
  $k$-elements for any $k<n$. 

Let $G$ be any $4$-generated group of exponent $p$ and class at
most~$2$, $p$ an odd prime. If $G$ is cyclic and nontrivial, then it
is neither capable nor extra-special. If $G$ is minimally $2$- or $3$-generated,
then it is capable. If $G$ is minimally $4$-generated, then it is
capable in all cases except when $G$ is extra-special of order $p^5$.
Therefore, we have
established the following:

\begin{theorem}
Let $G$ be a $4$-generated $p$-group of class at most two and
exponent~$p$. Then $G$ satisfies one and only one of the following:
\begin{itemize}
\item[(i)] $G$ is capable; or
\item[(ii)] $G$ is cyclic and nontrivial; or
\item[(iii)] $G$ is extra-special of order $p^5$ and exponent~$p$.
\end{itemize}
\label{th:classifnfour}
\end{theorem}

\begin{remark} It may be profitable to consider the following: each 
  $\langle\mathbf{u}\rangle \in Gr(1,U)$ determines, via $\langle
  \mathbf{u}\rangle^*$, an $(n-1)$-dimensional subspace of~$V$. For
  every $k\geq n-1$, this subspace determines a well-known subset (in
  fact, a subvariety) of $Gr(k,V)$, known as \textit{Schubert cells}
  or \textit{Schubert cycles}. For the case $n=4$, the above theorem
  states that the $X\in Gr(5,V)$ that are closed are exactly the ones
  that are in some Schubert cycle determined as
  $\langle\mathbf{u}\rangle$ ranges over all elements of $Gr(1,U)$.
\end{remark}

\section*{References}
\bibliographystyle{amsalpha}

\end{document}